\documentclass{article}
\usepackage{pinlabel}

\usepackage{graphicx}

\newcommand{\pic}[2]{\raisebox{-.5\height}{\includegraphics[scale=#2]{#1}}}
\newcommand{\pica}[2]{\raisebox{-.6\height}{\includegraphics[scale=#2]{#1}}}

\def\AB{\pic{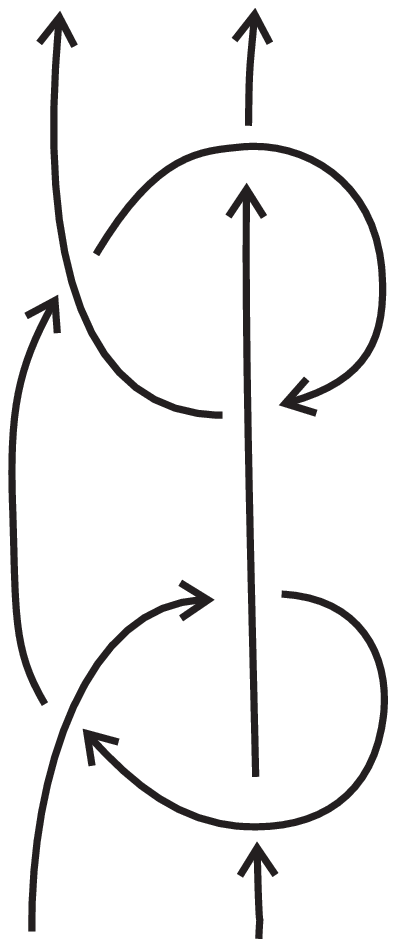}{.300}}
\def\A{\pic{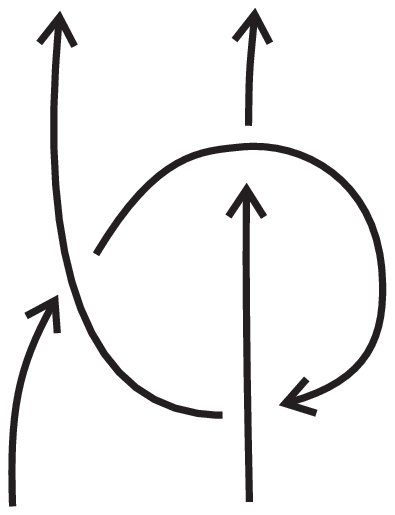}{.300}}
\def\B{\pic{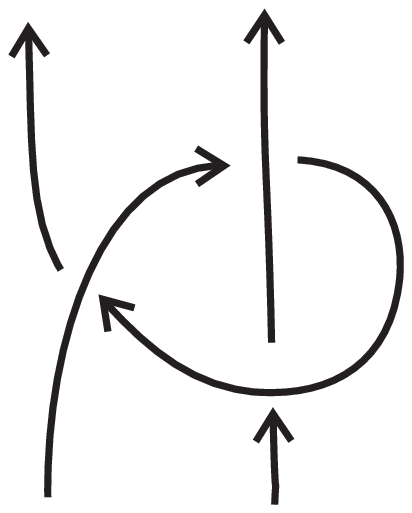}{.300}}
\def\ABreef{\pic{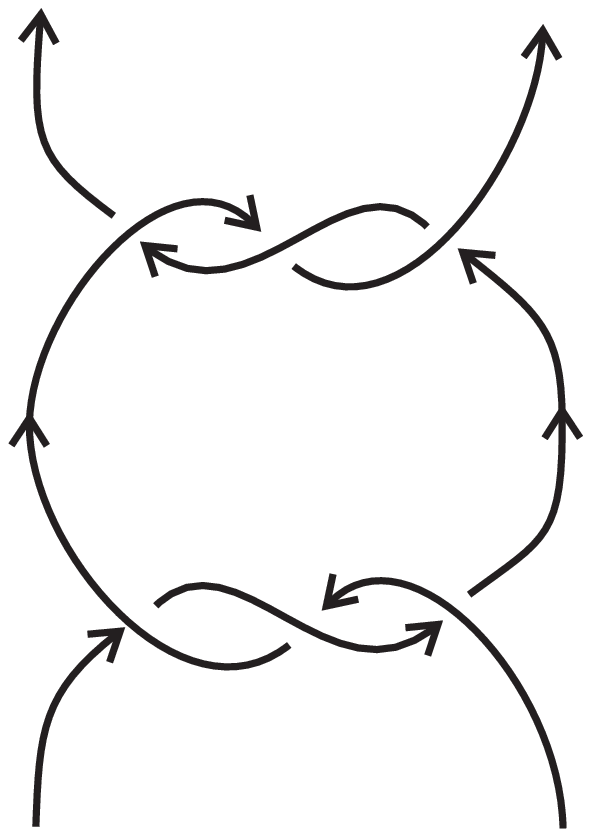}{.300}}

\def\closure{\pic{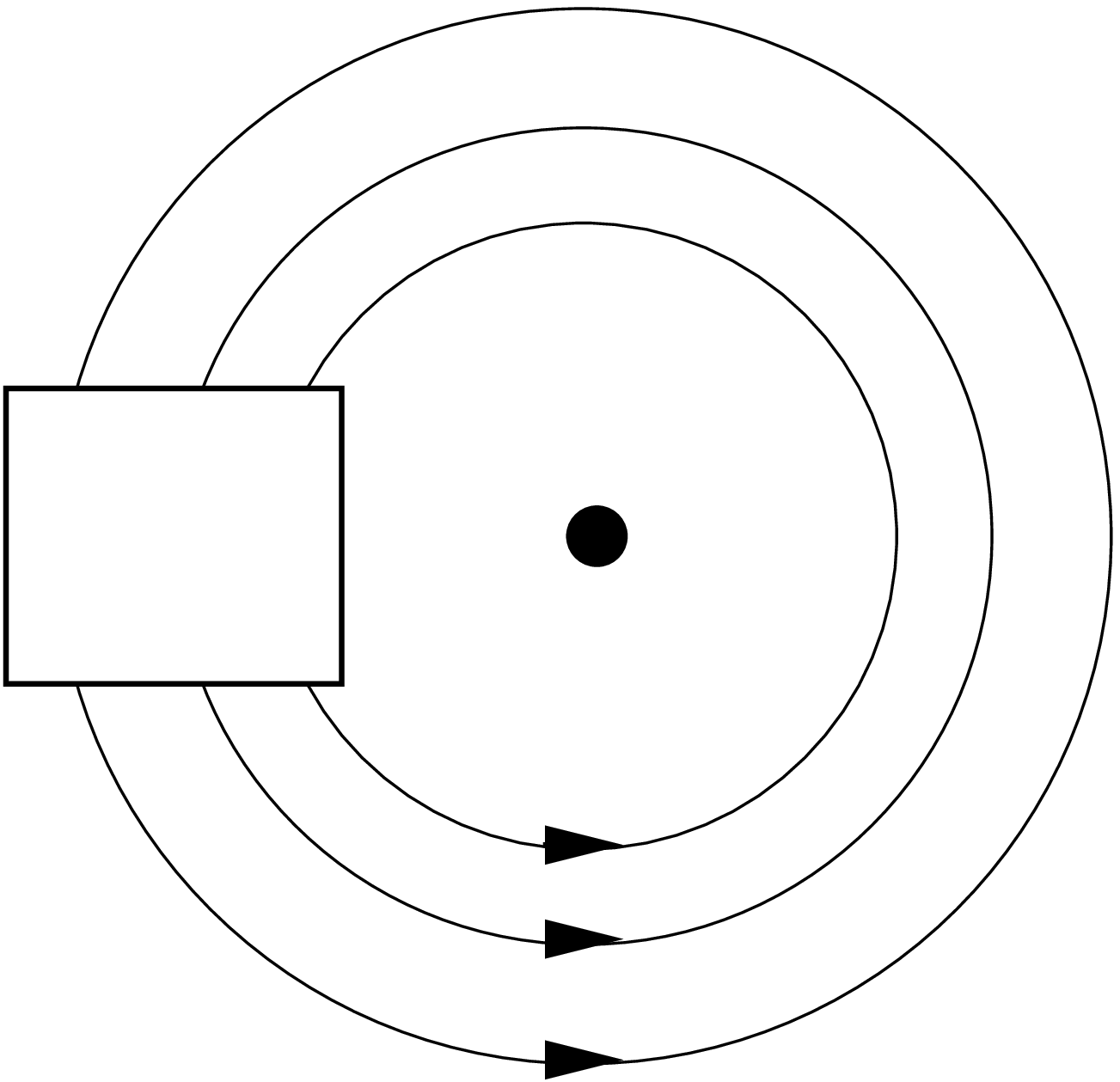}{.250}}

\def\Xor{\pic{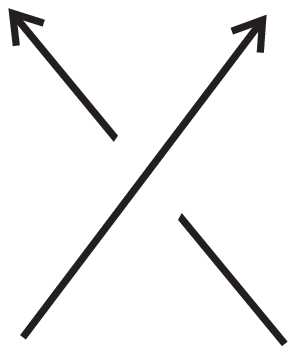}{.2500}}
\def\Yor{\pic{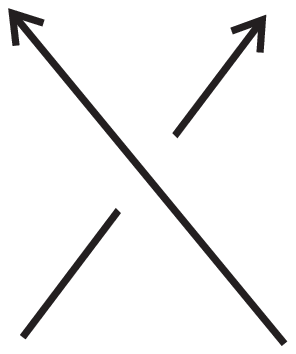} {.2500}}
\def\Ior{\pic{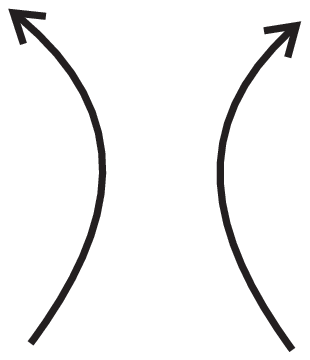} {.2500}}

\def\Satellite{\pic{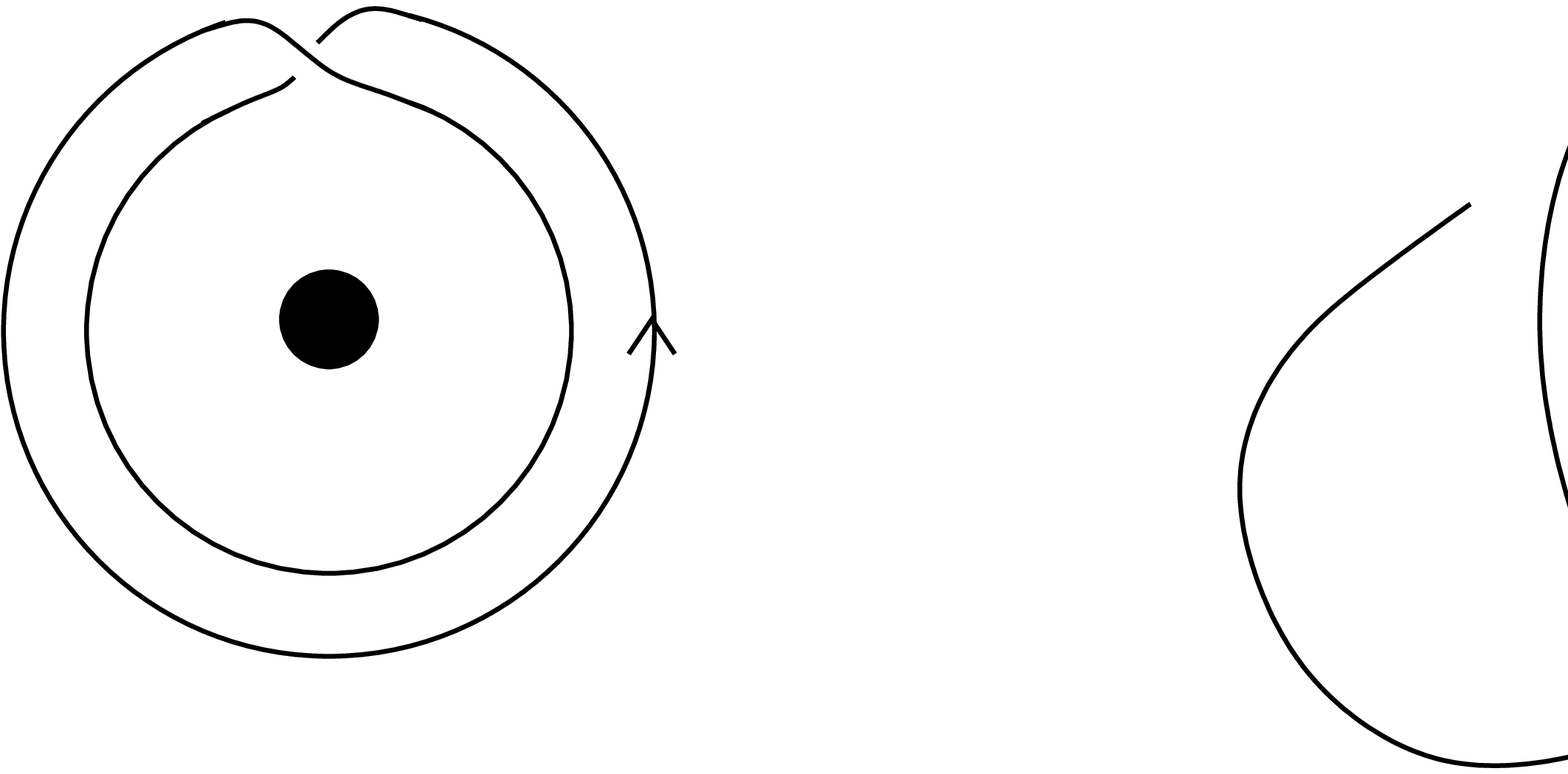}{.200}}

\def\Rmatrix{\pic{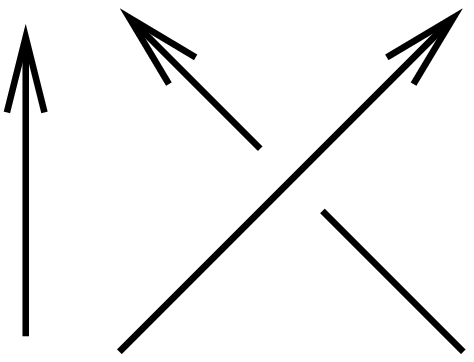} {.250}}
\def\threeparts{\pic{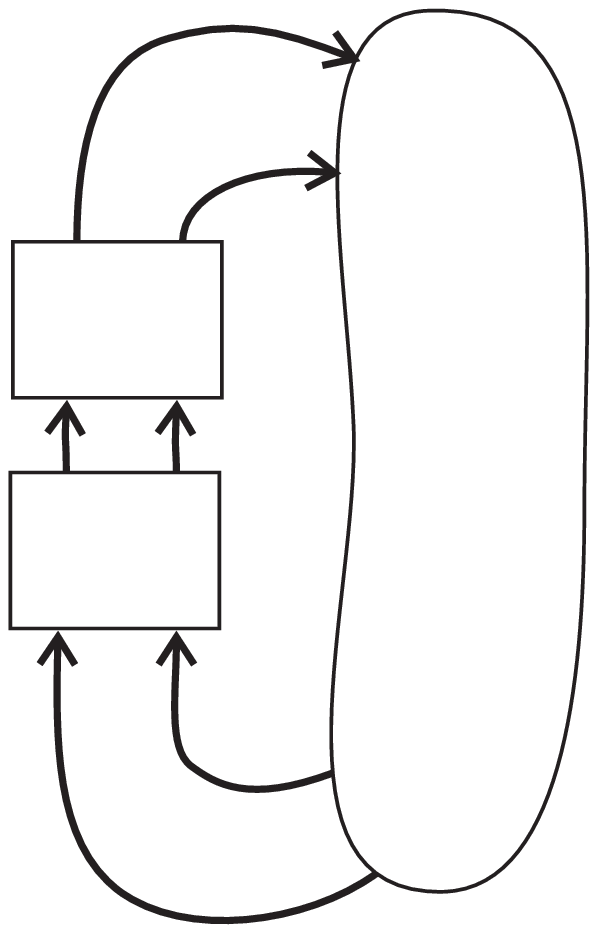}{.40}}
\def\unknot{\pic{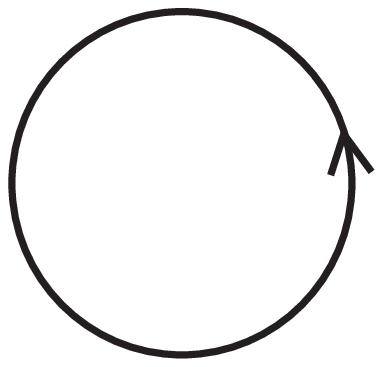} {.150}}
\def\Idor{\pic{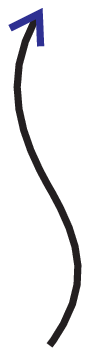} {.2500}}
\def\Rcurlor{\pic{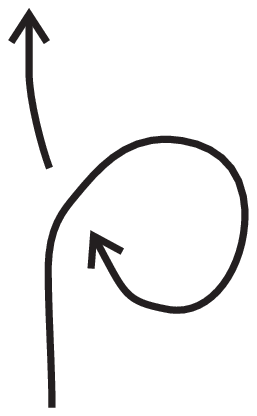} {.2500}}
\def\Conway{\pic{} {  .400}}
\def\KT{\pic{} {.400}}

\def\mutantone{\pica{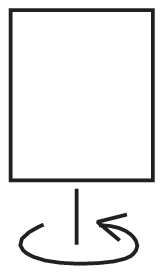} {.6}}
\def\mutanttwo{\pic{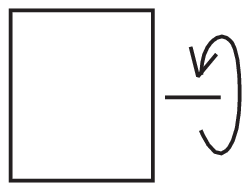} {.6}}
\def\mutantthree{\pic{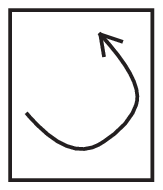} { .6}}

\def\mutantbox{\pic{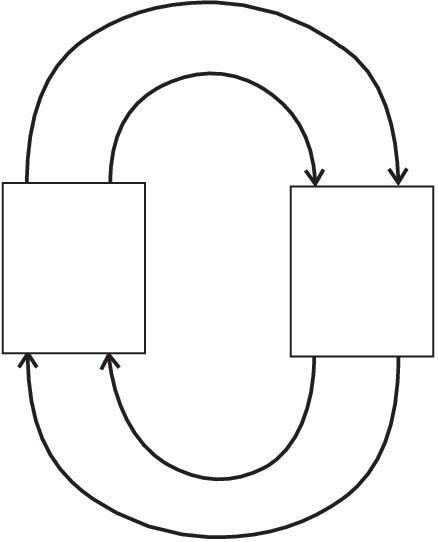}{ .6}}

\def\pretzel{\pic{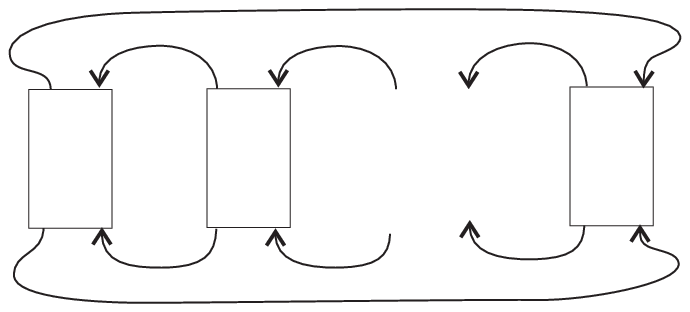}{.60}}
\def\pretzelgranny{\pic{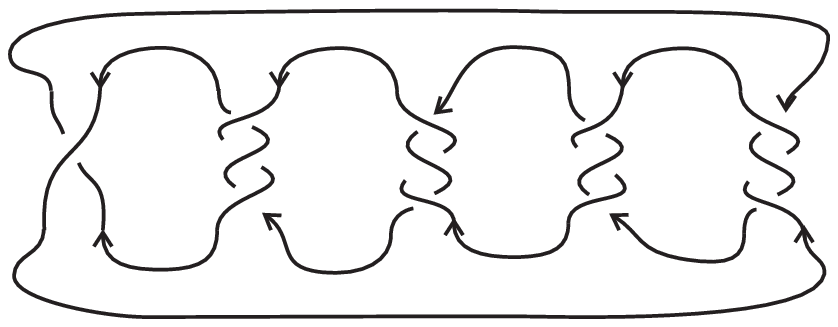}{.600}}
\def\pretzelreef{\pic{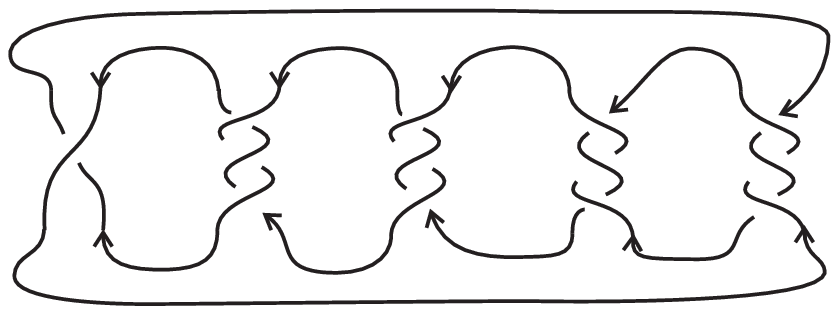}{.600}}

\usepackage{amssymb}
\usepackage{amsthm}

\newcommand{\bc}{\begin{center}}
\newcommand{\ec}{\end{center}}

\newcommand{\be}{\begin{equation}}
\newcommand{\ee}{\end{equation}}

\newtheorem{theorem}{Theorem} 
\newtheorem{corollary}[theorem]{Corollary}

\newtheorem{remark}{Remark}

\def\Oneone{\pic{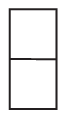} {.300}}

\def\Twoone{\pic{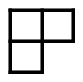} {.400}}
\def\Twotwo{\pic{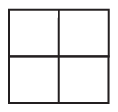}{.300}}
\def\Two{\pic{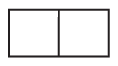}{.300}}
 
\def\Fourtwo{\pic{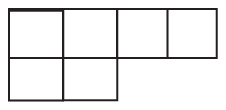}{.300}}

\def\Crossing{\pic{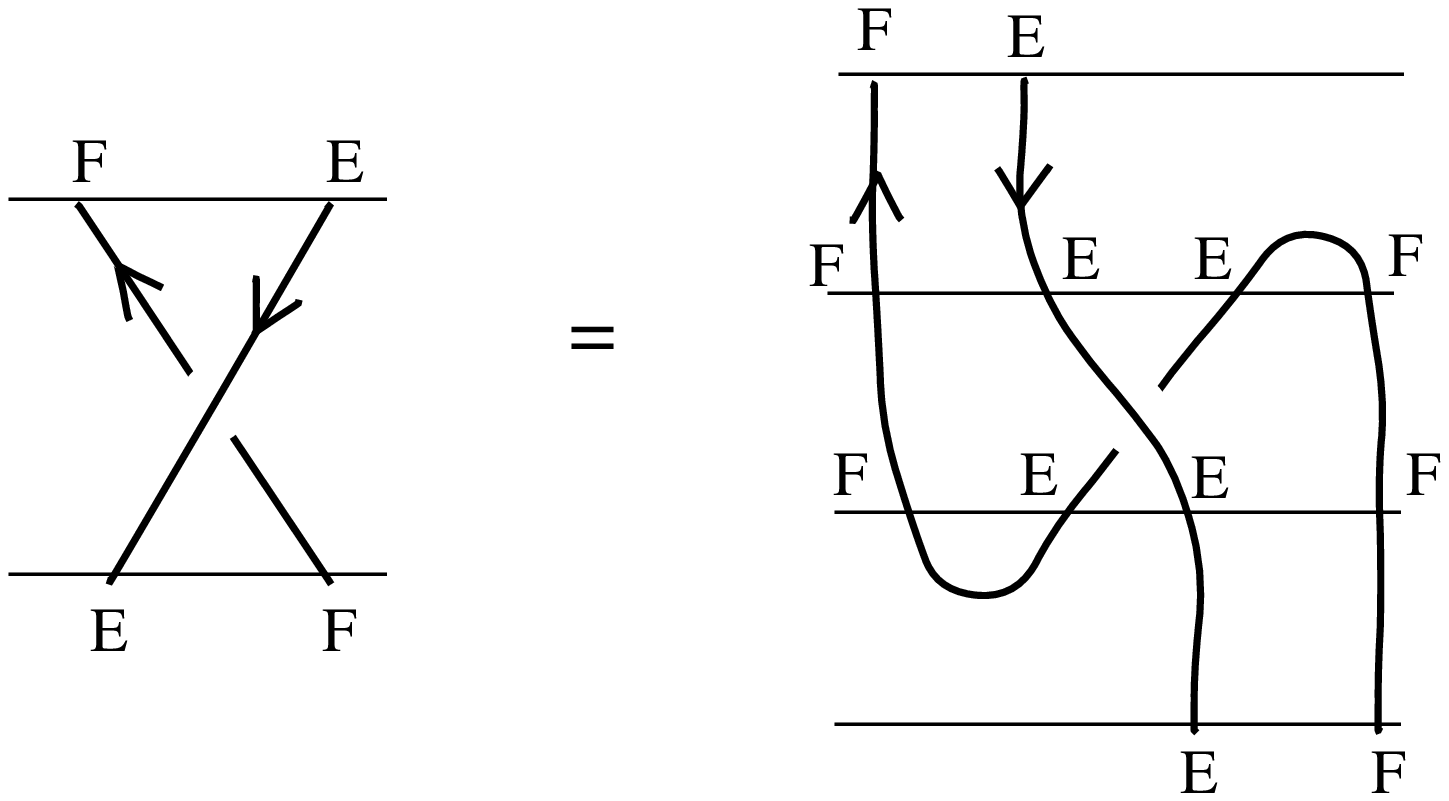} {.400}}

\def\cl{\centerline}
\def\gl{\lambda}
\def\x{\times}
\def\Oplus{\bigoplus}
\def\ds{\displaystyle}   
\def\tr{\makebox{tr}}
\begin{document}

\cl{\Large\bf Mutant knots with symmetry }                

\cl{\sc H.R.Morton}               
\vglue 0.1truein

{\small\sl 
\cl{Department of Mathematical Sciences,}           
\cl{University of Liverpool,}
\cl{Peach Street, Liverpool L69 7ZL, UK.}
}
%

\begin{abstract}
Mutant knots, in the sense of Conway, are known to share the same Homfly polynomial. Their $2$-string satellites also share the same Homfly polynomial, but in general their $m$-string satellites can have different Homfly polynomials for $m>2$. We show that, under conditions of extra symmetry on the constituent $2$-tangles, the directed $m$-string satellites of mutants share the same Homfly polynomial for $m<6$ in general, and for all choices of $m$ when the satellite is based on a cable knot pattern.

We give examples  of mutants with extra symmetry whose Homfly polynomials of some $6$-string satellites are different, by comparing their quantum $sl(3)$ invariants.
\end{abstract}

\section{Introduction}

This paper has been inspired by recent observations of Ochiai and Jun Murakami about the Homfly skein theory of $m$-parallels of certain symmetrical $2$-tangles. In \cite{Ochiai} Ochiai remarks that the $3$-parallels of the tangle $AB$ in figure \ref{figone} and its mirror image $\overline{A}\overline{B}=BA$ are equal in the Homfly skein of $6$-tangles, in other words, in the Hecke algebra $H_6$, \cite{AistonMorton}.
\begin{figure}[ht!]
{\labellist
\small
\pinlabel {$A$} at 156 647
\pinlabel {$B$} at 156 517
\endlabellist
\bc
\AB\qquad=\qquad \ABreef
\ec}
\caption{} \label{figone}
\end{figure}

As a consequence, the $3$-parallels of any mutant pair of knots given by composing the $2$-tangles $AB$ and $BA$ with any other $2$-tangle $C$ and then closing will share the same Homfly polynomial. 
  
This is in contrast with the known fact that $3$-parallels of mutant knots in general can have different Homfly polynomials, \cite{MortonTraczyk, MortonCromwell}.

There is interest in the extent to which the Homfly polynomial of $m$-parallels or other $m$-string satellites can distinguish mutants which are closures of $ABC$ and $BAC$ with $A$ and $B$ as above. Ochiai has found that the $4$-parallels of $AB$ and $BA$ are different in the skein $H_8$. 

\begin{figure}[ht]
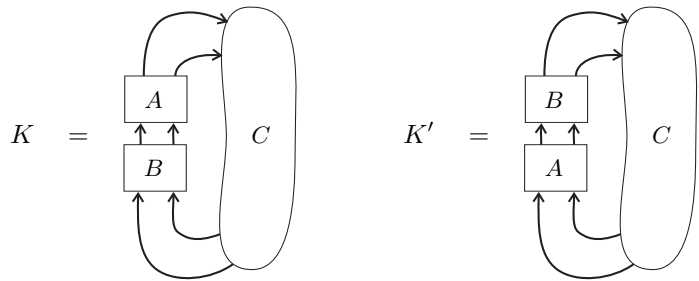

\bc 
$$ K \quad=\quad{
\labellist
\small
\pinlabel{{$A$}} at 226 509
\pinlabel{{$B$}} at  226 443
\pinlabel{{$C$}} at 330 475
\endlabellist
\threeparts
}
\qquad \qquad
K'\quad =\quad 
{
\labellist
\small
\pinlabel{{$B$}} at 226 509
\pinlabel{{$A$}} at  226 443
\pinlabel{{$C$}} at 330 475
\endlabellist
\threeparts
}
$$
\ec\caption{Tangle interchange}\label{threeparts}
\end{figure}

The purpose of this paper is to show that  if $A$ and $B$ are any two oriented $2$-tangles with symmetry
\bc
$ A$\quad = \quad
{\labellist
\pinlabel $A$ at 129 709

\endlabellist
\mutantone}\ \ , \quad $ B$\quad = \quad
{\labellist
\pinlabel $B$ at 129 709

\endlabellist
\mutantone}
\ec
  then the $m$-parallels, and indeed any directed $m$-string satellite, of knots $\widehat{\phantom A}ABC$ and $\widehat{\phantom A}BAC$ shown in figure \ref{threeparts} share the same Homfly polynomial for $m<6$. 

In contrast there exist examples of $A,B$ and $C$, including Ochiai's case with \bc $A\quad =$\quad \A \ ,\quad $B\quad =$\quad \B\ ,\ec
 for which the Homfly polynomials of the $6$-fold parallel are different.

As an unexpected extension  of the main result we show that the Homfly polynomial of a genuine connected cable, based on the $(m,n)$ torus knot pattern, with $m$ and $n$ coprime, for any number of strings, $m$, will not distinguish mutants with symmetry above, although a more general connected satellite pattern can do so.

The examples which exhibit differences for the directly oriented $6$-parallel can also be used to show that the $4$-parallels with two pairs of reverse strands have distinct Homfly polynomials.

The proofs are based on the relation of the Homfly satellite invariants to quantum $sl(N)$ invariants, and the techniques are an extension of work with Cromwell \cite{MortonCromwell} and with H. Ryder \cite{MortonRyder}. The eventual calculations that exhibit the difference of invariants in the specific example  depend on the $27$ dimensional irreducible module over $sl(3)$ corresponding to the partition $4,2$, and some Maple calculations following similar lines to those in \cite{MortonRyder}.

\section{Shared invariants of mutants}

 The term \emph{mutant} was coined by Conway, and  refers to the following
general construction.

Suppose that a knot $K$ can be decomposed into 
two oriented  $2$-tangles $F$ and $G$ 

\begin{center}
$K$ \quad {=} \quad
{
\labellist
\pinlabel{$F$} at 100 733
\pinlabel{$G$} at 181 733
\endlabellist
\mutantbox}
\end{center}

 A new knot $K'$
can be formed by replacing the  tangle $F$ with the tangle
$F'=\tau_i(F)$ given by rotating 
$F$ through $\pi $ in one of three ways,  

\bc
$\tau_1(F)$\quad {=} \quad {\labellist
\pinlabel{$F$} at 96 668
\endlabellist
\mutanttwo }\ ,\quad  $\tau_2(F)$\quad {=} \quad{\labellist
\pinlabel{$F$} at 102 629
\endlabellist
\mutantthree}\ ,\quad $\tau_3(F)$\quad {=} \quad
{\labellist
\pinlabel{$F$} at 129 709
\endlabellist
\mutantone }\ ,
\ec
reversing its string orientations if
necessary.  
Any of the three knots \bc 
$K'$\quad {=} \quad{\labellist
\pinlabel{{$\tau_i(F)$}} at 102 735
\pinlabel{{$G$}} at 181 735
\endlabellist
\mutantbox}
\ec
 is called a 
 \emph{mutant} of $K$. 

The two  $11$-crossing knots, $C$ and $KT$, with trivial  Alexander polynomial found by 
Conway and  Kinoshita-Teresaka are the best-known example of  mutant knots.

\medskip
\begin{center}{ $C\ =\ $ \Conway\  \qquad
 $KT\ =\ $\KT\ }
\end{center}
\subsection{Satellites}
\begin{figure}[ht]
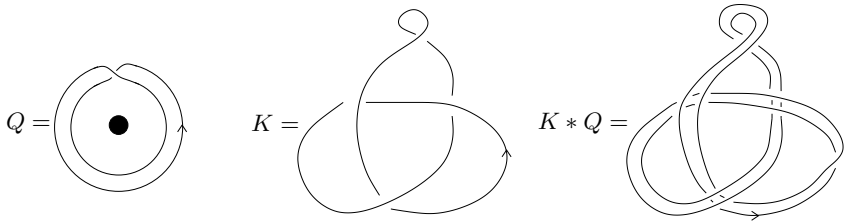

{\labellist
\small
\pinlabel {$Q=$} at -50 195
\pinlabel {$K=$} at 440 195
\pinlabel {$K*Q=$} at 1050 195
\endlabellist
\bc
\Satellite
\ec}
\caption{Satellite construction} \label{figsatellite}
\end{figure}

 A satellite of $K$ is determined by choosing a diagram $Q$ in the standard annulus, and then drawing $Q$ on the annular neighbourhood of $K$ determined by the framing, to give the satellite knot $K*Q$. We refer to this construction as {\em decorating $K$ with the pattern $Q$}, as shown in figure \ref{figsatellite}.

For fixed $Q$ the Homfly polynomial $P(K*Q)$ of the satellite is an invariant of the framed knot $K$. The invariants $P(K*Q)$ as $Q$ varies make up the   {\em Homfly satellite invariants} of $K$. We use the alternate notation $P(K;Q)$ in place of $P(K*Q)$ when we want to emphasise the dependence on $K$.

The general symmetry result compares the invariants of two knots $K$ and $K'$ made up of $2$-tangles $A$, $B$ and $C$,  by interchanging $A$ and $B$ as in figure \ref{threeparts}. 

\begin{theorem} \label{symmetry} Suppose that 
 $A$ and $B$ are both symmetric under the half-twist $\tau_3$, so that \bc
$ A$\quad = \quad
{\labellist
\pinlabel $A$ at 129 709

\endlabellist
\mutantone}\ \ , \quad $ B$\quad = \quad
{\labellist
\pinlabel $B$ at 129 709

\endlabellist
\mutantone}
\ec
Let $K$ and $K'$ be knots which are the closure of $ABC$ and $BAC$ respectively for any tangle $C$, as in figure \ref{threeparts}.
Then
$P(K*Q)=P(K'*Q)$ for every closed braid pattern $Q$ on $m<6$ strings.
\end{theorem}

\begin{remark} Our proof will apply equally to the case where $Q$ is the closure of a directly oriented $m$-tangle with $m<6$.
\end{remark}

In order to prove the theorem we must rewrite the Homfly satellite invariants in terms of quantum $sl(N)$ invariants, so we now give a brief summary of the relations bewteen these invariants, originally established by Wenzl. Further details can be found in \cite{AistonMorton} and the thesis of Lukac, \cite{Lukacthesis}, including details of  variant Homfly skeins with a framing correction factor, $x$. These are isomorphic to the skeins used here but the parameter allows a careful adjustment of the quadratic skein relation to agree directly with the natural relation arising from use of the quantum groups $sl(N)$.

\subsection{Homfly skeins}

For a surface $F$ with some designated input and output boundary points
the (linear) {Homfly skein} of $F$ is defined as linear combinations of oriented diagrams in $F$, up to Reidemeister moves II and III,  modulo the skein relations
\bc\begin{enumerate}
\item \qquad{$\Xor\  -\ \Yor \qquad =\qquad{(s-s^{-1})}\quad\ \Ior \ ,$}
\item \qquad
{$ \Rcurlor \qquad=\qquad {v^{-1}}\quad \Idor\ . $}
\end{enumerate}
\ec

It is an immediate consequence that 
\bc  \unknot  \ \Idor\quad=\quad $\delta$ \ \Idor,\ec
where $\delta =\ds\frac{v^{-1}-v}{s-s^{-1}}\in\Lambda$.
The coefficient ring $\Lambda$ is  taken as $Z[v^{\pm1},s^{\pm1}]$, with denominators $ s^r-s^{-r}, r\ge1$.

The skein of the annulus is denoted by ${\cal C}$. It becomes a commutative algebra with a product induced by placing one annulus outside another. 

The skein of the rectangle with $m$ inputs at the top and $m$ outputs at the bottom is denoted by $H_m$. We define a product in $H_m$ by stacking one rectangle above the other, obtaining the Hecke algebra $H_m(z)$, when $z=s-s^{-1}$ and the coefficients are extended to $\Lambda$. The Hecke algebra $H_m$ can   also be regarded as the group algebra of Artin's braid group $B_m$ generated by the  elementary braids $\sigma _i$, $i=1, \dots , m-1$, modulo the further quadratic relation $\sigma _i^2=z\sigma _i +1$.

\begin{figure}[ht!]
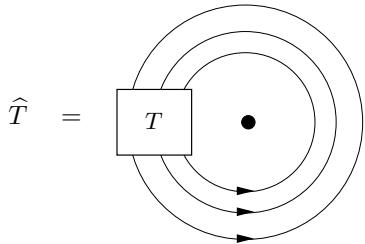

\labellist
\small
\pinlabel $T$ at 58 185
\endlabellist
  \begin{center}
$\widehat{T}\quad =\quad \closure $
  \end{center}
\caption{The closure map}\label{closuremap}
\end{figure}

The closure map from $H_m$ to ${\cal C}$ is the $\Lambda$-linear map induced by mapping a tangle $T$ to its closure $\widehat{T}$   in the annulus (see figure \ref{closuremap}).  We refer to a diagram $Q=\widehat{T}$ as a \emph{directly oriented} pattern.

The image of this map is denoted by ${\cal C}_m$, which
  has a useful interpretation as the space of symmetric polynomials of degree $m$ in variables $x_1, \ldots , x_N$ for large enough $N$.  Moreover, the submodule ${\cal C}_+ \subset {\cal C}$ spanned by the union $\cup _{m \geq 0}\, {\cal C}_m$ is a subalgebra of ${\cal C}$ isomorphic to the algebra of the symmetric functions.

\subsection{Quantum invariants}
A quantum group $\cal G$ is an algebra over a formal power series ring ${\bf Q}[[h]]$, typically a deformed version of a classical Lie algebra. We write $q=e^h, s=e^{h/2}$ when working in $sl(N)_q$.
A finite dimensional module over $\cal G$ is a linear space on which $\cal G$ acts.

Crucially, $\cal G$ has a coproduct $\Delta$ which ensures that the tensor product $V\otimes W$ of two modules is also a module.
It also has a {\em universal $R$-matrix} (in a completion of ${\cal G}\otimes{\cal G}$) which determines a well-behaved module isomorphism $$R_{VW}:V\otimes W \to W\otimes V.$$

This has a diagrammatic view
indicating its use in converting coloured tangles to module homomorphisms.
\bc
{\labellist
\small
\pinlabel{{$ W\  \otimes \ V$}} at 84 156
\pinlabel{{$V\ \otimes \ W$}} at 84 12
\pinlabel{$R_{VW}$} at -50 84
\endlabellist}
\Rmatrix
\ec

A braid $\beta$ on $m$ strings with permutation $\pi\in S_m$ and a colouring of the strings by modules $V_1,\ldots,V_m$ leads to a module homomorphism 
$$J_\beta:V_1\otimes\cdots\otimes V_m \to V_{\pi(1)}\otimes\cdots\otimes V_{\pi(m)}$$ using $R_{V_i,V_j}^{\pm1}$ at each elementary braid crossing.
The homomorphism $J_\beta$ depends {\em only on the braid} $\beta$ itself, not its decomposition into crossings, by the Yang-Baxter relation for the universal $R$-matrix.

When $V_i=V$ for all $i$ we get a module homomorphism $J_\beta:W\to W$, where 
$W=V^{\otimes m}$. Equally, a directed $m$-tangle $T$ determines an endomorphism $J_T$ of $W=V^{\otimes m}$. Now any $sl(N)$ module $W$ decomposes as a  direct sum $\bigoplus {(W_\mu\otimes V_\mu^{(N)})}$, where $W_\mu$ is the  linear subspace consisting of the {\em highest weight vectors} of type $\mu$  associated to the module $V_\mu^{(N)}$.   Highest weight subspaces of each type are preserved by module homomorphisms, and so $J_T$ determines (and is determined by) the restrictions $J_T(\mu):W_\mu \to W_\mu$ for each $\mu$.

If a knot $K$ is decorated by a pattern $Q$ which is the closure of an $m$-tangle $T$ then its quantum invariant $J(K*Q;V)$ can be found from the endomorphism $J_T$ of $W=V^{\otimes m}$ in terms of the quantum invariants of $K$ and the highest weight maps $J_T(\mu):W_\mu \to W_\mu$   by the formula
\be J(K*Q;V)=\sum c_\mu J(K;V_\mu^{(N)}) \label{weighttrace}\ee with $c_\mu=\mbox{tr}\, J_T(\mu)$. This formula follows from lemma II.4.4 in Turaev's book \cite{Turaevbook}. Here $\mu $ runs over partitions with at most $N$ parts when we are working with $sl(N)$, and we set $c_\mu=0$ when $W$ has no highest weight vectors of type $\mu$.

\begin{proof}[Proof of theorem \ref{symmetry}]
Take $V=V^{(N)}$ as the fundamental module of dimension $N$ for $sl(N)$. Then the only highest weight types $\mu$ which occur in equation (\ref{weighttrace}) are partitions of $m$ with at most $N$ rows. Because $J(K*Q;V^{(N)})=P(K*Q)$ when $v=s^{-N}$ we can show that $P(K*Q)=P(K'*Q)$ by showing that $J(K*Q;V^{(N)})=J(K'*Q;V^{(N)})$ for all $N$. By equation \ref{weighttrace} it is then enough to show that $J(K;V_\mu^{(N)})= J(K';V_\mu^{(N)})$ for all $N$ and all partitions $\mu\vdash m$.

Now each tangle $A$ and $B$ determines an endomorphism $J_A, J_B$ of $V_\mu\otimes V_\mu$. If $J_A$ and $J_B$ commute then $J(K;V_\mu)=J(K';V_\mu)$. The endomorphisms $J_A$ and $J_B$ are determined by their restriction $J_A(\nu), J_B(\nu)$ to the highest weight subspaces $W_\nu$ in the decomposition $V_\mu\otimes V_\mu=\sum W_\nu\otimes V_\nu$, so it is enough to show that $J_A(\nu)$ and $J_B(\nu)$ commute where $V_\nu$ is a summand of $V_\mu\otimes V_\mu$. This is certainly the case for all $\nu$ where $W_\nu$ is $1$-dimensional, which includes the case of single row or column partitions $\mu$, \cite{MortonCromwell}.

As a special case of the work of Rosso and Jones, \cite{Rosso,MortonManchon}, we know that the endomorphism  of $V_\mu\otimes V_\mu$ for the full twist $\Delta^2$ on two strings operates as a scalar $e^{f(\nu)}$ on each highest weight space $W_\nu$, while the half twist $\Delta$, represented by the $R$-matrix $R_{V_\mu V_\mu}$, operates on $W_\nu$ with two eigenvalues $\pm e^{\frac{1}{2}f(\nu)}$.

The positive and negative eigenspaces corrspond to the classical decomposition of the Schur function $(s_\mu)^2$ into symmetric and skew-symmetric parts, $h_2(s_\mu)$ and $e_2(s_\mu)$, and the dimension of each eigenspace of $W_\nu$ is the multiplicity of $s_\nu$ in $h_2(s_\mu)$ and $e_2(s_\mu)$ respectively.

Now $A=\tau_3(A)$, so that $A\Delta=\Delta A$. Hence the endomorphism $J_A$, and similarly $J_B$, preserves the positive and negative eigenspaces of each $W_\nu$. If these eigenspaces   have dimension $1$ or $0$ then $J_A$ and $J_B$ will commute on $W_\nu$.

The theorem is then established by checking that no $s_\nu$ occurs in $h_2(s_\mu)$ or $e_2(s_\mu)$ with multiplicity $>1$ for any $\mu$ with $|\mu|\le 5$. The decomposition of all of these can be quickly confirmed using the Maple program SF of Stembridge \cite{Stembridge}.
\end{proof}
\begin{corollary}

Examples  include $k$-pretzel knots {$K(a_1,\ldots,a_k)$} with odd $a_i$. 

\bc
{\labellist
\small 
\pinlabel{{$a_1$}} at 143 670
\pinlabel{{$a_2$}} at 194 672
\pinlabel{{$a_k$}} at 300 674
\endlabellist
\pretzel}
\ec

Here the numbers $a_i$ can be permuted without changing the Homfly polynomial of any satellite with $\le5$-strings.

\end{corollary}

\section{Satellites with different Homfly polynomials}

A further check with the program SF when $|\mu|=6$ shows that there are just three partitions, $\mu=4,2$, its conjugate $\mu = 2,2,1,1$ and $\mu=3,2,1$ whose symmetric square $h_2[s_\mu]$ contains summands with multiplicity $>1$, as does the exterior squares of $\mu=3,2,1$. Explicitly $h_2[s_{4,\,2}] = {s_{8, \,4}} + {s_{8, \,2, \,2}} + {s_{7, \,4, \,1}} + {s_{7, \,3
, \,2}} + {s_{7, \,3, \,1, \,1}} + {s_{6, \,6}} + {s_{6, \,5, \,1
}} + 2\,{s_{6, \,4, \,2}} + {s_{6, \,3, \,2, \,1}} + {s_{6, \,2, 
\,2, \,2}} + {s_{5, \,5, \,1, \,1}} 
 + {s_{5, \,4, \,3}} + {s_{5, \,4, \,2, \,1}} + {s_{5, \,3
, \,3, \,1}} + {s_{4, \,4, \,4}} + {s_{4, \,4, \,2, \,2}}$.  This means that, although $m$-string satellites of $K$ and $K'$ must share the Homfly polynomial when $m\le 5$, it is possible for the Homfly polynomials of some $6$-string satellites to differ.

We give an example now where this does indeed happen.

\begin{theorem}

Let $K$ and $K'$ be the pretzel knots $K=K(1,3,3,-3,-3)$ and $K'=K(1,3,-3,3,-3)$. 

\bc \pretzelreef\qquad\pretzelgranny
\ec

The $6$-fold parallels $K*Q$ and $K'*Q$, where $Q$ is the closure of the identity braid on $6$ strings, have different Homfly polynomials.
\end{theorem}

\begin{proof}Write $K$ and $K'$ as the closure of the products $\Delta ABAB$ and $\Delta BAAB$ respectively, where \bc $A\quad =$\quad \A \ ,\quad $B\quad =$\quad \B\ ,\ec
 are  the  partially closed $3$-braids shown, and $\Delta$ is the positive half-twist. We show that $P(K*Q)\ne P(K'*Q)$ when $v=s^{-3}$. These values are given by the $sl(3)$ quantum invariants $J(K*Q;V^{(3)})$ and $J(K'*Q;V^{(3)})$, where $V^{(3)}$ is the fundamental $3$-dimensional module   for $sl(3)$.  Since $Q$ is the closure of the identity braid on $6$ strings it induces the identity endomorphism on the module $(V^{(3)})^{\otimes 6}$. This module decomposes as $\Oplus W_\mu\otimes V_\mu^{(3)}$ where $\mu$ runs through partitions of $6 $ with at most $3$ rows. The trace of the identity on $W_\mu$ is just $d_\mu=\dim W_\mu$, giving $$J(K*Q;V^{(3)})=\sum d_\mu J(K;V_\mu^{(3)}).$$
 
The only partition $\mu$ in this range for which the exterior or symmetric square contains highest weight vectors of multiplicity $>1$ is the partition $\mu=4,2$, since the partition $\mu=2,2,1,1$ has $4$ rows and the repeated factors for $\mu=3,2,1$ occur for partitions with more than $3$ rows.
Now $J_A(\mu)J_B(\mu)=J_B(\mu)J_A(\mu)$ for all other $\mu$ since $A$ and $B$ are symmetric up to altering the framing on both strings, while maintaining the writhe. Then
$$P(K*Q)- P(K'*Q)=d_\mu(J(K;V_\mu^{(3)})-J(K';V_\mu^{(3)}))$$ when $v=s^{-3}$ and $\mu=4,2$. Since $d_\mu\ne0$ it is enough to show that $J(K;V_\mu^{(3)})\ne J(K';V_\mu^{(3)})$.
 The  module $V_\mu^{(3)}$ has dimension $27$.

We   now work in the quantum group $sl(3)$ and  drop the superscript $(3)$ from the irreducible modules.

Decompose the module $V_\mu\otimes V_\mu$ as $\sum W_\nu\otimes V_\nu$ and compare the endomorphisms given by the tangles $T=ABAB\Delta$ and $T'=BAAB\Delta$.

In this case just one of the invariant subspaces of highest weight vectors has dimension $> 1$. It can be shown that the corresponding $2\x2$ matrices $A_\mu$ and $B_\mu$ arising from the two mirror-image tangles $A$ and $B$ with $3$ crossings satisfy $\tr (A_\mu B_\mu A_\mu B_\mu -A_\mu A_\mu B_\mu B_\mu)\ne 0$, which results in a difference in their $sl(3)$ invariants $J(K;V_\gl)$. 

 None of the other $6$-cell invariants differ on the two knots. Consequently the $6$-parallels have different $sl(3)$ invariants. The $sl(3)$ invariant of the $6$-parallels of the two pretzel knots coloured with the fundamental module, and thus  their Homfly polynomials, are then different.
\end{proof}

\subsection{Use of the quantum group $sl(3)_q$}

The calculation of the $2\x2$ matrices $A_\nu$ and $B_\nu$ giving the effect of the two tangles on the highest weight vectors where there is a 2-dimensional highest weight subspace of the symmetric part of the module depends on finding the explicit action of the quantum group on the 27-dimensional module $V_{\mu}^{(3)}$ with $\mu=4,2$ and its tensor square, as well as the homomorphism representing its $R$-matrix.  I used the linear algebra packages in Maple to handle the matrix working and subsequent polynomial factorisation, following fairly closely the techniques developed with H. Ryder in the paper \cite{MortonRyder}. 

In the interests of reproducibility I give an account of the methods used, and some of the checks applied during the calculations, to test against known properties.

We start from a presentation of the quantum group $sl(3)_q$ as an algebra with
six generators, $X_1^{\pm},\,X_2^{\pm},\,H_1,\,H_2$, and a description of the
comultiplication and antipode. 

Let $M$ be any finite-dimensional left module
over $sl(3)_q$. The action of any one of these six generators $Y$  will
determine a linear endomorphism $Y_M$ of $M$. We  build up explicit matrices
for these endomorphisms on a selection of low-dimensional modules, using the
comultiplication to deal with the tensor product of two known modules, and the
antipode to construct the action on the linear dual of a known module. We must
eventually determine the matrices $Y_M$ for our module 
$M=V_{\Fourtwo}$, and find the $729\times729$  $R$-matrix, $R_{MM}$ 
which represents the endomorphism of $M\otimes M$ needed for crossings.

We follow Kassel in the basic description of the quantum group from
using generators $H_1$ and $H_2$ for the Cartan
sub-algebra, but with generators $X_i^\pm$ in place of $X_i$ and $Y_i$. We  use
the notation  $K_i=\exp(hH_i/4)$, and  set $a=\exp(h/4), \, s=\exp(h/2)=a^2$
and $q=\exp(h)=s^2$, unlike Kassel. The generators satisfy the commutation
relations 
$$[H_i,H_j]=0, \ [H_i,X_j^\pm]=\pm a_{ij}X_j^\pm, \ 
[X_i^+,X_i^-]=(K_i^2-K_i^{-2})/(s-s^{-1}),$$ 
where
$(a_{ij})=\pmatrix{2&-1\cr-1&2\cr}$ is the Cartan matrix for $SU(3)$ (and
also the Serre relations of degree 3 between $X_1^\pm$ and $X_2^\pm$).

Comultiplication is given by
\[\begin{array}{rl}\Delta(H_i)&=H_i\otimes I+I\otimes H_i,\\ (\hbox{so
}\Delta(K_i)&= K_i\otimes K_i,)\\
\Delta(X_i^\pm)&=X_i^\pm\otimes K_i+ K_i^{-1}\otimes X_i^\pm,\\
\end{array}
\] and the antipode $S$ by $S(X_i^\pm) =-s^{\pm 1}X_i^\pm$, $S(H_i)=-H_i$, $
S(K_i)=K_i^{-1}$.

The fundamental $3$-dimensional module, which we denote by $E$, has a basis in
which the quantum group generators are represented by the matrices $Y_E$ as
listed here.
\[X_1^+=\pmatrix{0&1&0\cr0&0&0\cr0&0&0\cr},\
X_2^+=\pmatrix{0&0&0\cr0&0&1\cr0&0&0\cr}\]
\[ X_1^-=\pmatrix{0&0&0\cr1&0&0\cr0&0&0\cr},\
X_2^-=\pmatrix{0&0&0\cr0&0&0\cr0&1&0\cr}\]
\[H_1=\pmatrix{1&0&0\cr0&-1&0\cr0&0&0\cr},\
=\pmatrix{0&0&0\cr0&1&0\cr0&0&-1\cr}.\]

For calculations we   keep track of the elements $K_i$ rather than $H_i$,
represented by
\[K_1=\pmatrix{a&0&0\cr0&a^{-1}&0\cr0&0&1\cr},\
K_2=\pmatrix{1&0&0\cr0&a&0\cr0&0&a^{-1}\cr}\] for the module $E$.

We can then write down the elements $Y_{EE}$ for the actions of the generators
$Y$ on the module $E\otimes E$, from the comultiplication formulae. The
$R$-matrix $R_{EE}$  can be given, up to a scalar,
by the prescription
\[\begin{array}{rl}R_{EE}(e_i\otimes e_j)&=e_j\otimes e_i, \hbox{ if }i>j,\\
&=s\,e_i\otimes e_i, \hbox{ if } i=j,\\ &=e_j\otimes e_i+(s-s^{-1})e_i\otimes
e_j, \hbox{ if }i<j,
\\
\end{array}
\] for basis elements $\{e_i\}$ of $E$.

The linear dual $M^*$ of a module $M$ becomes a module when the action of a
generator $Y$ on $f\in M^*$ is defined by $<Y_{M^*}f,v>=<f,S(Y_M)v>$, for $v\in
M$. For the dual module $F=E^*$ we then have matrices for $Y_F$, relative to
the dual basis, as follows.

\[X_1^+=\pmatrix{0&0&0\cr-s&0&0\cr0&0&0\cr},\
X_2^+=\pmatrix{0&0&0\cr0&0&0\cr0&-s&0\cr}\]
\[ X_1^-=\pmatrix{0&-s^{-1}&0\cr0&0&0\cr0&0&0\cr},\
X_2^-=\pmatrix{0&0&0\cr0&0&-s^{-1}\cr0&0&0\cr}\]
\[K_1=\pmatrix{a^{-1}&0&0\cr0&a&0\cr0&0&1\cr},\
K_2=\pmatrix{1&0&0\cr0&a^{-1}&0\cr0&0&a\cr}.\]

The most reliable way to work out the $R$-matrices $R_{EF}, R_{FE}$ and
$R_{FF}$ is to combine $R_{EE}$ with module homomorphisms $\hbox{cup}_{EF}$, $
\hbox{cup}_{FE}$, $\hbox{cap}_{EF}$ and $\hbox{cap}_{FE}$ between the modules
$E\otimes F$, $F\otimes E$ and the trivial 1-dimensional module, $I$, on which
$X_i^\pm$ acts as zero and $K_i$ as the identity. The matrices are determined up to a
scalar by such considerations; a choice for one dictates the rest.

Once these matrices have been found they can be combined with the matrix
$R_{EE}^{-1}$ to construct the $R$-matrices $R_{EF},R_{FE},R_{FF}$, using the
diagram shown below, for example, to determine $R_{EF}$. This gives 
\[ R_{EF}=(1_F\otimes 1_E\otimes  \hbox{cap}_{EF})\circ 
(1_F\otimes R_{EE}^{-1}\otimes 1_F)\circ (\hbox{cup}_{FE}\otimes
1_E\otimes 1_F) .\]
\begin{center}
\Crossing
\end{center}

 The module structure
of 
$M=V_{\Fourtwo}$ can  be found by identifying $M$ as a
$27$-dimensional submodule of  $V_{\Twotwo}\otimes V_{\Two}$, while the two $6$-dimensional modules $V_{\Two}$ and $V_{\Twotwo}$ are themselves submodules of $E\otimes E$ and $F\otimes F$ respectively.

We know, by the Pieri formula, that there is
a direct sum decomposition of $V_{\Twotwo}\otimes V_{\Two}$ as $M\oplus N$, where $M=V_{\Fourtwo}$ and  $N$ is the sum of the $8$-dimensional module $V_{\Twoone}$ and the $1$-dimensional trivial module.

We first identify  the module $V_{\Two}$ as a submodule of $E\otimes E$, knowing that $E\otimes E$ is isomorphic  to $V_{\Two} \otimes F$.
 The full twist element on the two strings both coloured by $E$ is represented by $R^2_{EE}$ which acts on $E\otimes E$ as a scalar on each of the two irreducible submodules $V_{\Two}$ and $F$.

Use Maple to find bases for the two eigenspaces of $R^2_{EE}$. Then we can identify $V_{\Two}$ with the $6$-dimensional one, and write $P$ and $Q$ for the $9\x6$ and $9\x3$ matrices whose columns are these bases. The partitioned matrix $(P|Q)$
is invertible, and its inverse, found by Maple, can be written as
$\ds\left(\ds{R\over S}\right)$, where  $R$ is a $6\x9$ matrix
with $RP=I_{6}$ and $RQ=0$.

Regard $P=\mbox{inj}M_1EE$ as the matrix representing the
inclusion of the module $V_{\Two}$ into $E\otimes E$. Then  $R=\mbox{proj}EEM_1$ is the
matrix, in the same basis, of the projection from $E\otimes E$ to $V_{\Two}$.
For $M_1 =V_{\Two}$ the module generators $Y_{M_1}$   are given by  $Y_{M_1}=R\,Y_{EE}\,P$, giving the explicit
action of the quantum group on $V_{\Two}$.  

We perform a similar calculation on $F\otimes F$ to identify the module $M_2=V_{\Twotwo}$ and the matrices $\mbox{inj}M_2FF$ and $\mbox{proj}FFM_2$, giving the action of the quantum group on $M_2=V_{\Twotwo}$  in a similar way.

We use   inclusion and projection further to find the four $6^2\times 6^2$
$R$-matrices $R_{M_iM_j}$.  For example, to construct $R_{M_1M_2}:M_1\otimes M_2\to M_2\otimes M_1$, first map $M_1\otimes M_2$ to $E\otimes E\otimes F\otimes F$ by $\mbox{inj}M_1EE\otimes \mbox{inj}M_2FF$. Then construct the $R$-matrix crossing two strings with $E\otimes E$ and two with $F\otimes F$ as the composite of $1\otimes R_{EF}\otimes 1$ , $R_{EF}\otimes R_{FE}$ and $1\otimes R_{FF}\otimes 1$, and finally compose with the projections  $\mbox{proj}FFM_2 \otimes \mbox{proj}EEM_1$.

A similar calculation on the module $M_1\otimes M_2$ yields the submodule $M=V_{\Fourtwo}$. The full twist on two strings, one coloured by $M_1$ and one by $M_2$, is represented by the product $R_{M_2M_1}R_{M_1M_2}$ and will have one 27-dimensional eigenspace $M$ complemented by two other eigenspaces. Taking the bases of these eigenspaces in a partitioned $36\x36 $ matrix as above will determine a $36\x 27$ matrix $P=\mbox{inj}MM_1M_2$ and a $27\x36$ matrix $R=\mbox{proj}M_1M_2M$. The quantum group actions $Y_{M_1M_2}$ on the tensor product are determined by the coproduct formulae, and the actions $Y_M$ are then given from these using $P$ and $R$. These in turn give rise to the quantum group actions $Y_{MM}$ on $M\otimes M$.

We are also able to construct the $27^2\x27^2$  $R$-matrix $R_{MM}$ using the same inclusion and projection to map $M\otimes M$ into $M_1\otimes M_2\otimes M_1\otimes M_2$, followed by the matrix for crossing four strands, built up from the $R$-matrices $R_{M_iM_j}$ and then the projections back to $M\otimes M$.

\subsection{Completing the calculations}\label{calculation}

\begin{remark}
We can reach this stage directly if we know the six module generators $Y_M$ and the $R$-matrix $R_{MM}$ for the module $M=V_{\Fourtwo}$. We can then calculate the module generators $Y_{MM}$ using the coproduct, and the twisting element $T_M=(K_{1M})^4(K_{2M})^4$.
\end{remark}

Knowing the module generators  $Y_{MM}$ gives an immediate means of finding the highest weight vectors as common null-vectors of $X^+_{iMM}$, and their weights can be identified.  All the submodules of $M\otimes M$ occur with multiplicity $1$ except $V_\nu$ with partition $\nu=6,4,2$ whose highest weights are $2,2$.
The $3$-dimensional space $W_\nu$ of highest weight vectors for $\nu$ is found by solving the linear equations $X^+_{1MM}v=0$, $X^+_{2MM}v=0$, $K_{1MM}v=a^2v$ and $K_{2MM}v=a^2v$  for $v$. We then find the $2$-dimensional positive eigenspace for $R_{MM}$ on $W_\nu$. The endomorphisms $J_A$ and $J_B$ will preserve this eigenspace.

  Represent the
$3$-braid $\sigma_2 \sigma_1^{-1}\sigma_2$ in the $2$-tangle $A$ by an endomorphism $F_A$ of $M\otimes M\otimes M$,
using $R_{MM}$ and its inverse. Then use $T_M$ and the partial trace to close
off one string, hence giving the endomorphism $J_A$ of $M\otimes M$
determined by $A$. Explicitly, choose a basis $\{e_i\}$ of $M$  and write $$F_A(v\otimes T_M(e_i))=\sum_j f_{ij}(v)\otimes e_j$$ with $f_{ij}(v)\in M\otimes M$. Then $J_A(v)=\sum_i f_{ii}(v)$. Applied to each of the two vectors in the highest weight space this   determines a $2\x2$ matrix $A_\nu$ representing the restriction of $J_A$ to this subspace. Similarly $B_\nu$ is found using the mirror image braid $\sigma_2^{-1} \sigma_1\sigma_2^{-1}$.

We know that $R_{MM}$ acts as a scalar on the $2$-dimensional space  so 
$J(K;V_\mu)-J(K';V_\mu)$ is a non-zero scalar multiple of $\tr(A_\nu B_\nu A_\nu B_\nu-B_\nu A_\nu A_\nu B_\nu)$. 

This difference is $2(q^6+q^5+q^4+q^3+q^2+q+1)(q^4+1)(q^6+q^3+1)^2(q^4-q^2+1)^2(q^4
+q^3+q^2+q+1)^3(q^2+1)^4(q^2+q+1)^4(q^2-q+1)^4(q+1)^{10}(q-1)^{18}$, up to a power of $q=s^2$ and the quantum dimension of $V_\nu$.

\subsection{Further examples of difference}

Using the same matrices $A_\nu$ and $B_\nu$ it is possible to find further pretzel knot examples based on sequences of the tangles $A$ and $B$ where the $6$-parallels have different Homfly polynomial, such as the knots $K(3,3,3,-3,-3)$ and $K(3,3,-3,3,-3)$. The difference here is the same as for the first example multiplied by the factor $2q^{32}-q^{31}-3q^{30}+5q^{29}+3q^{28}-10q^{27}+q^{26}+14q^{25}-6q^{24}-19q^{23}+
21q^{22}+20q^{21}-46q^{20}+2q^{19}+61q^{18}-48q^{17}-35q^{16}+83q^{15}-27q^{14}
-66q^{13}+72q^{12}+3q^{11}-57q^{10}+40q^9+10q^8-33q^7+16q^6+7q^5-12q
^4+7q^3-4q+2$. The same calculations   guarantee that satellites based on any  closed $6$-tangle $Q=\widehat T$ will have different Homfly polynomial, provided that the trace $c_\mu$ of the endomorphism $J_{\widehat T} $ on the highest weight space $W_\mu$ of $V^{\otimes 6}$ is non-zero, where $\mu$ is the partition $4,2$. This will be the case for most, but not all, patterns $Q$, and certainly will be the case for many satellites which are knots rather than links.

The calculations in section \ref{calculation} also show that the $4$-parallels of the two pretzel knots $K(1,3,3,-3,-3)$ and $K(1,3,-3,3,-3)$ with two strings oriented in one direction and two in the opposite direction will have different Homfly polynomials, by using the decomposition of the corresponding $sl(3)_q$ module $W=V\otimes V\otimes V_{\Oneone}\otimes V_{\Oneone}$ into a sum of irreducible $sl(3)_q$ modules. The only module to figure in this decomposition with any multiplicity in its symmetric or exterior square is again $V_{\Fourtwo}$. The calculations above, using the fact that Homfly with $v=s^{-3}$ can be calculated by colouring strings with reverse orientation by the dual module $V^*$ to the fundamental module, and that this is $V_{\Oneone}$ for $sl(3)_q$.

\section{Cable patterns}
By way of contrast, if the pattern $Q$ is a cable on any number of strings then $K*Q$ and $K'*Q$ share the same Homfly polynomial, where $K$ and $K'$ have the same symmetry as in theorem \ref{symmetry}.

\begin{theorem} \label{cablesymmetry} Suppose that 
 $A$ and $B$ are both symmetric under the half-twist $\tau_3$, so that \bc
$ A$\quad = \quad
{\labellist
\pinlabel $A$ at 129 709

\endlabellist
\mutantone}\ \ , \quad $ B$\quad = \quad
{\labellist
\pinlabel $B$ at 129 709

\endlabellist
\mutantone}
\ec
Let $K$ and $K'$ be knots which are the closure of $ABC$ and $BAC$ respectively for any tangle $C$, as in figure \ref{threeparts}.
Then
$P(K*Q)=P(K'*Q)$ for every $(m,n)$ cable pattern $Q$ where $m$ and $n$ are coprime. \end{theorem}
\begin{proof} As in the proof of theorem \ref{symmetry} we show that $J(K*Q;V^{(N)})=J(K'*Q;V^{(N)})$ for all $N$. By equation \ref{weighttrace} it is then enough to show that $J(K;V_\mu^{(N)})= J(K';V_\mu^{(N)})$ for all $N$ and all partitions $\mu\vdash m$ for which the coefficient $c_\mu\ne 0$. The coefficients $c_\mu$ depend on the pattern $Q$ and arise as the trace of the endomorphism $J_T$ when restricted to the highest weight space $W_\mu \subset V^{\otimes m}$, where $Q$ is the closure of the $m$-braid $T=(\sigma_1\sigma_2\cdots\sigma_{m-1})^n$.  

It is shown in \cite{Rosso}, (see also \cite{MortonManchon}), that for any such   cable $Q$ the only non-zero coefficients $c_\mu$ occur when the partition $\mu$ is  a {\em hook}, if  $m$ and $n$ are coprime .
It is then enough to show that $J(K;V_\mu^{(N)})= J(K';V_\mu^{(N)})$ for all hook partitions $\mu$.

Using the same argument as in theorem \ref{symmetry} it remains to check that 
  no Schur function $s_\nu$ occurs with multiplicity $>1$ in the decomposition of either the symmetric or exterior squares, $h_2(s_\mu)$ or $e_2(s_\mu)$, for any hook partition $\mu$. This fact has been established by Carbonara, Remmel and Yang in  theorem 3 of  \cite{Carbonara}, and so the proof is complete.
\end{proof}
\begin{remark}
Theorem \ref{cablesymmetry} highlights the importance of a precise terminology for different types of satellite. The term {\em cable} is sometimes used to mean any satellite, while there is a clear distiction here between the behaviour of cables and of parallels or other satellites, which is not primarily a matter of the number of components of the satellite.
\end{remark}

\section*{Acknowledgements}

I would like to thank the Topology group at Universidad Complutense, Madrid, for their hospitality during some of the preparation of this article.

I am grateful to Bernard Leclerc and Jean-Yves Thibon for help in identifying the decomposition of the symmetric and antisymmetric square of Schur functions of degree $\le6$, which figure in the arguments. I would also 
like to thank  John Stembridge for further suggestions of methods for establishing the general result about cables, and Christine Bessenrodt for help in tracking down the article \cite{Carbonara} used in proving theorem \ref{cablesymmetry}.

\noindent May 2007


\begin{thebibliography}{99}

\bibitem{AistonMorton} AK Aiston and HR Morton. Idempotents of Hecke algebras of type $A$. \emph {J. Knot Theory  Ramifications} {\textbf  7} (1998), 463--487.

\bibitem{Carbonara} JO Carbonara, JB Remmel and  M Yang. Exact formulas for the plethysm $s_2 [s _{(1 ^a; b)}]$ and $s_{1^ 2} [s_{(1^a; b)}]$. Technical report, Mathematical Sciences Institute, Cornell University, 1992.

\bibitem{Lukacthesis} SG Lukac. Homfly skeins and the Hopf link.  PhD. thesis, University of Liverpool, 2001.

\bibitem{MortonCromwell} HR Morton and PR Cromwell.
Distinguishing mutants by knot polynomials.
\emph{ J. Knot Theory Ramifications} {\textbf  5} (1996), 225--238. 

\bibitem{MortonManchon} HR Morton and PMG Manchon. Some basic formulas in the Homfly skein of the annulus. Preprint, University of Liverpool 2007.

\bibitem{MortonRyder} HR Morton and HJ Ryder:
Mutants and $SU(3)q$ invariants. 
In `Geometry and Topology Monographs', Vol.1: The Epstein Birthday Schrift. (1998), 365--381. 


\bibitem{MortonTraczyk} HR Morton and P Traczyk.
The Jones polynomial of satellite links around mutants. 
In `Braids', ed. Joan S. Birman and Anatoly Libgober, \emph{Contemporary Mathematics} \textbf {78}, Amer. Math. Soc. (1988), 587--592. 


 \bibitem{Ochiai} M Ochiai and N Morimura. Base tangle decompositions of $n$-string tangles with $1<n<10$. Preprint 2006.


\bibitem{Rosso} M Rosso and VFR Jones. {On the invariants of torus knots derived from quantum groups. } \emph {J. Knot Theory  Ramifications } {\textbf 2}  (1993), 97--112.

\bibitem{Stembridge} JR Stembridge. A Maple package for symmetric functions. Version 2.4, (2005), University of Michigan, {\tt www.math.lsa.umich.edu/\~ { }jrs}

\bibitem{Turaevbook} VG Turaev. { Quantum invariants of knots and 3-manifolds.}
De Gruyter Studies in Mathematics, 18.
Walter de Gruyter and Co., Berlin, 1994.

\end{thebibliography}
\end{document}